\documentclass[%
  twocolumn
   , colorlinks 
]{mpi2015-cscpreprint}

\usepackage[american]{babel}

\usepackage{graphicx}

\usepackage{amssymb}
\usepackage{amsthm}

\usepackage{wrapfig}
\usepackage{tikz}
\usepackage{pgfplots}
\pgfplotsset{compat=1.18}
\usepackage{filecontents}
\usepackage{subfigure}

\usepackage{wrapfig}
\usepackage{cuted}
\usepackage{hyperref}
\usepackage[linesnumbered,ruled,vlined]{algorithm2e}

\AtBeginDocument{

}

\begin{document}


\title{SIGGRAPH: G: Improved Projective Dynamics Global Using Snapshots-based Reduced Bases}

\author[$\ast$]{Shaimaa Monem}
\affil[$\ast$]{Max Planck Institute for Dynamics of Complex Technical Systems, Magdeburg, Germany.\authorcr%
  \email{monem@mpi-magdeburg.mpg.de}, \orcid{0009-0008-4038-3452}}

\author[$\ast$]{Peter Benner}
\affil[$\ast$]{Max Planck Institute for Dynamics of Complex Technical Systems, Magdeburg, Germany.\authorcr%
	\email{benner@mpi-magdeburg.mpg.de}, \orcid{0000-0003-3362-4103}}

\author[$\dagger$]{Christian Lessig}
\affil[$\dagger$]{European Centre for Medium-Range Weather Forecasts, Bonn, Germany.\authorcr%
	\email{christian.lessig@ovgu.de}, \orcid{0000-0002-2740-6815}}

\shorttitle{IPDGS}
\shortauthor{Monem et. al.}
\shortdate{}



\abstract{%

  We propose a snapshots-based method to compute reduction subspaces for physics-based simulations. Our method is applicable to any mesh with some artistic prior knowledge of the solution and only requires a record of existing solutions during, for instance, the range-of-motion test that is required before approving a mesh character for an application. Our subspaces span a wider range of motion, especially large deformations, and rotations by default. Compared to the state-of-the-art, we achieve improved numerical stability, computational efficiency,  and more realistic simulations with a smaller sub-space.}

\novelty{We introduce a novel technique to compute reduced subspace that has the capability to preserve global rotation information, for real-time interactive physics-based simulations of deformable objects for applications in computer graphics.}

\maketitle

\begin{figure*}[t!]
	\centering
	\includegraphics[width=\textwidth, height=0.5\textwidth]{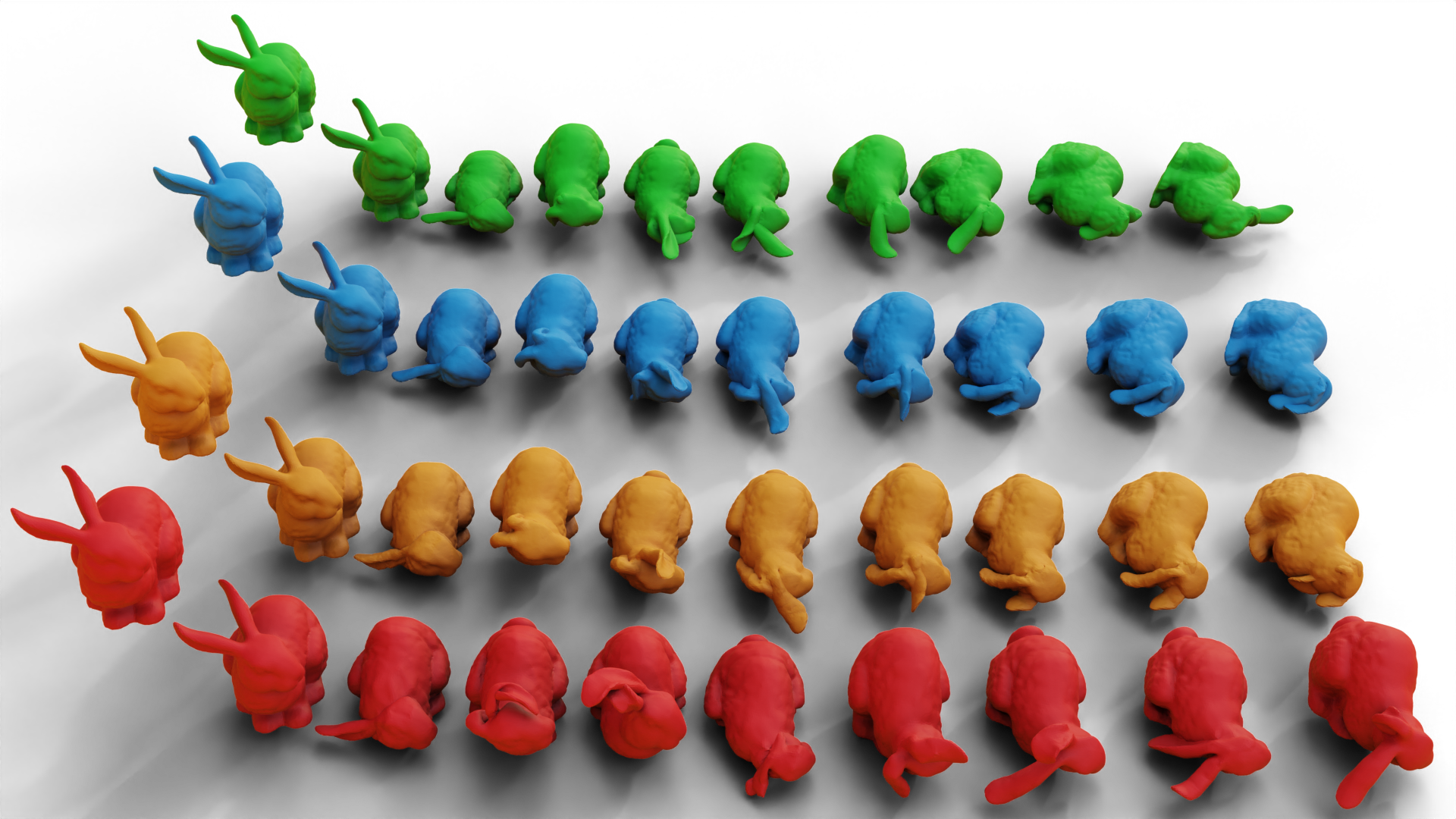}
	\caption{Graphical comparison showcase of full simulations (displayed in green) alongside with reduced models, localized PCA (in blue), SPLOCS (in orange), and LBS (in red). The series includes centered frames from a simulation rendering a "falling bunny."}
	\label{fig:FULL_PCA_SPLOCS_LBS_posOnly}
\end{figure*}

\section{Introduction}%
\label{sec:intro}
	
	Model reduction techniques for real-time physics-based simulations of elastic bodies are well studied in the literature, with a wide range of applications in engineering \cite{BennerBook22}, fluid dynamics \cite{lassila14}, computer graphics \cite{Trusty23}, the game industry \cite{Bnneland19}, bio-systems \cite{Radulescu12}, and more. Deformable objects present considerable computational challenges because of their complex nonlinear behavior and the high dimensionality of the problems they pose. However, solid deformables have the benefit of requiring fewer frequencies to accurately represent their range of deformation compared to fluids, which need higher levels of detail and hence larger subspaces.
	
	A variety of methods have been developed to enhance the speed of physics simulations while still preserving their accuracy and stability. These methods encompass explicit and implicit finite element methods \cite{snowden2017methods}, position-based dynamics \cite{PBS15, PBSSurvay14}, and projective dynamics \cite{Bouaziz14}. Each technique brings unique benefits and drawbacks in terms of speed, precision, and robustness, providing various compromises that can meet the specific needs of different applications.
	
	This work is inspired by, first, hyper-reduced projective dynamics \cite{Brandt18} which introduced linear blend skinning reduced subspaces, and achieved notable speed improvements while retaining high accuracy compared to other methods, and second, by sparse localized components \cite{Neumann13}.
	Another promising approach for model reduction is based on snapshot methods, such as sparse proper orthogonal bases and sparse localized components, and other snapshots methods which pinpoint the most crucial deformation modes from a collection of snapshots gathered from pre-computed simulations or experiments \cite{Barbic05, Stanton13}. In our study, we started by observing the missing factors in the currently available reduction techniques:
	
	\begin{itemize}
		\item Reduction subspaces are mainly designed based on mass and stiffness system matrices, as well as on rest state. Although they provide essential information about the simulation dynamics, they mostly lack information about large deformation that might be encountered during simulations, especially when user real-time interaction is involved,
		\item Existing methods are very limited in their ability to express global and local rotations accurately.
	\end{itemize}
	
	All of the above-mentioned issues are present in the current state-of-the-art introduced by Brandt et. al. \cite{Brandt18} as linear blend skinning subspace, which is designed from rest state and skinning weights, see Figure (\ref{fig:FULL_PCA_SPLOCS_LBS_posOnly}). We discuss those factors further in Section (\ref{related_work}) in detail. In the current work, we slightly extend our previous focus on accelerating the evaluation of positions state \cite{Monem23}, and we apply our method to projective dynamics \cite{Bouaziz14} to compare to the current state of art \cite{Brandt18}. It is worth noting that snapshot methods can be adapted to any physics simulations as long as some information about the solution is available.

\section{Related Work}
\label{related_work}

	One effective technique to accelerate the simulation of deformable elastic bodies is by determining a small set of bases vectors that spans a subspace with much lower degrees of freedom that contains an approximate solution of the full high dimensional simulations under study. Speed up is gained when computations are projected onto and conducted within the lower dimensional subspace \cite{huang2019survey}. Studies to determine such bases are common in engineering, graphics, and computer animation  \cite{Pentland_Williams89, He01book, Barbivc11, Barbic05}. A big challenge in graphics is to compute a set of bases that incorporate both translation and rotation, and we show in our study that the snapshots method is worthwhile in terms of accuracy and computational costs.
	
	Brandt et al. \cite{Brandt18} employed linear blend skinning to create these reduced subspaces for positional and constraint data, achieving faster but less accurate simulations. Our research focuses on the reduction of vertex positions through various bases, evaluating these against skinning-derived subspaces. Moreover, we also can go further and sparsify as well as localize the bases in our approach, similar to sparse optimized components of Neumann et. al. \cite{Neumann13}, while we test both global and local bases.
	
	\subsection{Modal Analysis Bases}
	
	The pioneering work of \cite{Terzopoulos87} led to a clear theory and model for elastic deformable materials. Pentland and Williams \cite{Pentland_Williams89} introduced good vibration and reduction subspaces to computer graphics and animations applications, Choi and Ko \cite{Min_Ko05} showed that modal analysis based on the linear strain tensor is not suitable for rotation capture in case large deformation is involved. We demonstrate that rotations can be easily incorporated in the subspace with much less effort and steps.
	
	When modal analysis was first introduced by Pentland and Williams \cite{Pentland_Williams89}, modal bases were computed using mass and stiffness matrices, which solely depend on the initial state. Therefore, such bases do not carry much information about larger deformation that is not reflected in the rest state and so is the time integration which might suffer from shearing or warping artifacts and the absence of rotational information due to restrictive linearization. To compensate for the above-mentioned artifacts, different studies have been carried out. For instance, \cite{Capell02} exploits a volumetric control skeleton to divide the character mesh into overlapping regions associated with different bones, then runs locally linearized simulations and at the end hierarchically blends them before visualization. Similar to \cite{Mueller02}, \cite{Min_Ko05} handles rotation separately by keeping track of local rotations of all the vertices in the mesh, through warping stiffness matrix in global and local respectively. It is possible to extend modal bases to account for non-linearity, for example, \cite{Barbic05} incorporates tangent linear vibration modes by computing the mass-normalized derivatives of the Hessian stiffness tensor and blend them with the mass-orthonormal principal components, while \cite{vonTycowicz13} applies linear deformation to modal bases to extract additional modes that represent rotational factor and as well as geometrical and material properties.
	
	\subsection{Sparse and Local Bases}
	
	\paragraph{Linear Blend Skinning Bases} For rigged characters, reduction subspaces can also be constructed using linear blend skinning scalar weights and rest states \cite{Hahn12, Jacobson12}, where the solver computes the associated translations and rotations for every frame \cite{Brandt18, Brandt19}. Such reduction technique can be used to produce complementary dynamics when combined with skinning eigenmodes \cite{Benchekroun23}, which have shown to be adequate for heterogeneous elastodynamics simulations \cite{Trusty23}. The big advantage of these subspaces is that they are sparse and local by construction, which leads to a large reduction in computational costs.
	
	\paragraph{Localized Bases}
	Global bases computed through principal component analysis (PCA) \cite{Alexa00}, also known as proper orthogonal decomposition (POD) \cite{Hinze05} is well known to span be the best approximation for a full solution at any given rank. However, including global information can easily grow the computational burden. \cite{Neumann13} addressed this challenge, by enforcing sparsity and locality through optimization of the PCA bases computed on the sequence of animated meshes. The SPLOCS bases in \cite{Neumann13} not only tackle the computational complexity but they have meaningful interpretation as well, and can be used for 3D editing and shape deformation. 
	
	
	\section{Approach}
	A finite element model for an elastic object under forces, is governed by Newton´s second law of motion $$ \mathbf{M} \ddot{q} + \mathbf{D} \dot{q} + \mathbf{K} q = F,$$
	where $\mathbf{M}$, $\mathbf{D}$, $\mathbf{K}$ and $\mathbf{F}$ are the mass, damping, stiffness, and force matrices respectively and $q(t) \in \mathbb{R}^{n \times 3}$ represents vertex position displacements over time.
	
	Modal bases, the traditional PCA bases, also called displacement modes \cite{Pentland_Williams89}, were computed through the generalized eigenvalue problem $\mathbf{K} \Phi = \mathbf{M} \Phi \Lambda$, and only bases associated with low frequencies or high eigenvalues in $\Lambda$. The problem with this approach is that $\mathbf{M}$ and $\mathbf{K}$ are completely characterized by the rest state, and hence are the bases $\Phi$. Frame-dependent artifacts that were observed in previous work might aggravate over time. Furthermore, the reduced simulations cannot express large deformations that did not appear in the rest states. Hence they need to be combined with other techniques to bring back the missing rotational factor \cite{Mueller02, Min_Ko05, Capell02}.
	
	We also note that subspace reduction techniques for elastic objects that include snapshots start by eliminating global rotations with respect to average character mesh. This is done in the pre-processing steps of sparse PCA, sparse localized components (SPLOCS) \cite{Neumann13}, and co-rotational techniques \cite{huang2019survey}. We propose a novel approach that incorporates the rotational components and avoids the pre-alignment step, allowing us to construct bases or components suitable for real-time interactive applications that account for both translation and rotations simultaneously. Our approach achieves faster global linear solve, thanks to the sparsity of the bases, while maintaining better accuracy and stability compared to HRPD \cite{Brandt18}. 

\section{Method and Implementations}  

	The trajectory of vertex positions $q(t) \in \mathbb{R}^{n \times 3}$, without damping effects, is influenced by both internal forces $f_{in}(q) = \sum_j W_j(q)$ and external forces $f_{ext}$, 
	
	\begin{equation}
		M \ddot{q} = f_{in}(q) + f_{ext}.
		\label{eqn:Newtons}
	\end{equation}
	
	Exploiting the variational formulation of this mechanical system, the task becomes to solve a minimization problem of the form
	
	\begin{equation}\label{opt_eq}
		\underset{q_{k+1}}{\text{min}} \ \  \frac{1}{(\delta t)^2} \left( \| M^{\frac{1}{2}} ( q_{k+1} - s(q_{k}) )\|^2_F + \sum_j W_j(q_{k}) \right)
	\end{equation}
	
	to determine the next position $ q_{k+1}$. The optimization problem (\ref{opt_eq}) can be reformulated as a linear system solve. The projective dynamics method \cite{Bouaziz14} uses an ADMM solver that alternates between constraint projections $p(q_{k})$, through local solves, and the evaluation of new positions $q_{k+1}$ for the next frame, by solving global linear system at the fixed constraint projection, as in equation (\ref{eqn: trajectory}).   
	
	\begin{equation}
		\underbrace{\left( \frac{M}{(\delta t)^2} + \sum_j \omega_j S^T_j S_j \right)}_{A} q_{k+1} = \underbrace{ \overbrace{\frac{M \ s_k(q_k)}{(\delta t)^2}}^{{\color{red} B u(t)}} + \overbrace{\sum_j S^T_j p_j(q_k)}^{{\color{red} f(t)}} }_{b_k} 
		\label{eqn: trajectory}
	\end{equation}
	$p_j(q)$ represents the mapping of positions onto the $j^{th}$ constraint manifold $C_j$. Each constraint is weighted by a non-negative $w_j$, influenced by specific elements selected by matrix $S_j$, and calculated over a time increment $\delta t$. The overarching linear equation (\ref{eqn:trajectory}) can be expressed as $A q_{k+1} = B u(t) + f(t)$, where $A$ and $B$, calculated as $\frac{M}{(\delta t)^t}$, are the system matrices. Here, $u(t) = s_k(q_k)$ denotes the time-variant input, and $f(t)$ is a nonlinear function updated at every time step.
	
	For vertex position minimization, we aim to find a linear subspace $U \in \mathbb{R}^{n \times k}$ (with $k$ significantly smaller than $n$), which approximates $q$ as $U \tilde{q}$, where $\tilde{q} \in \mathbb{R}^{k \times 3}$. This leads to a more computationally feasible linear system, as we can solve for $\tilde{q}$ instead of $q$ which bypasses the need for a costly computation of equation (\ref{eqn: trajectory}) i. e. 
	
	\begin{equation}
		U^T \ E \ U \tilde{q}(t) =  U^T \ B u(t) + U^T \ f(U \ \tilde{q}(t)).
		\label{eqn: reducedSys}
	\end{equation}
	
\begin{figure}[h]
		\centering
		\vspace{-0.2cm}
		\includegraphics[trim={4.5cm 10cm 8.2cm 10cm}, clip, width=0.44\textwidth, keepaspectratio]{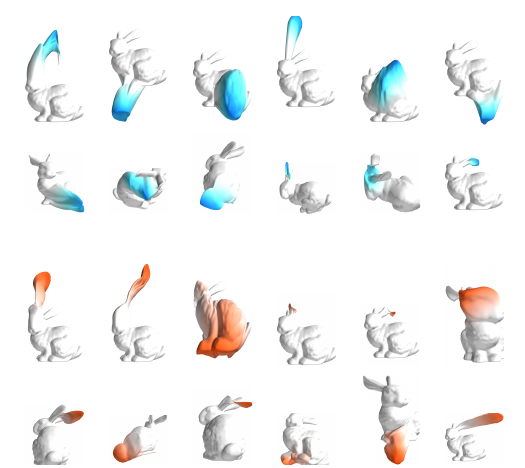}
		\caption{Few local components for a bunny deformed under gravitational forces: POD (blue) and SPLOCS (orange).}
		\label{fig:bunny_components_PCA_SPLOCS}
		
	\end{figure}

	Bases $U$ can be computed purely from skinning weights and rest position as in the current state-of-the-art \cite{Brandt18}. Running the implementations of \cite{Brandt18} shows that this vertex positions reduction method is highly unstable and comes with very obvious visual artifacts, see Figure (\ref{fig:bunny_bases_comparision}).
	We were inspired by techniques like PCA performed on entire simulation frames \cite{Alexa00, Neumann13}, which proved beneficial for reducing both the dimensions and the approximation error of simulated data. We computed weighted PCA bases and sparse localized bases for a simple experiment simulating the bunny mesh falling under gravitational forces, see Figure (\ref{fig:bunny_components_PCA_SPLOCS}). 
	
	Localized PCA bases show better numerical stability and efficiency for presenting large deformations, Figure (\ref{fig:bunny_bases_comparision}) illustrates this observation on the mesh layer, visual artifacts that can usually be covered by texture and not observed from a far camera perspective as in Figure (\ref{fig:FULL_PCA_SPLOCS_LBS_posOnly}). PCA bases can be computed with available snapshots, typically involving the capture of keyframes during motion tests—an essential step in validating character meshes for use in animations and interactive games. These frames effectively represent the artist's vision. Additionally, localized PCA bases can accurately depict motions not explicitly captured in the initial snapshots \cite{Neumann13}, eliminating the need for an extensive collection of data.
	
	In Algorithm (\ref{alg: PCA_SPLOCS}), the computation of PCA bases $U_{splocs}$ precedes the creation of SPLOCS bases. The global sparsification process via the ADMM solver, however, tends to remove the rotational elements of the bases, causing the mesh to appear fragmented and static, particularly in areas not covered by SPLOCS components, leading to noticeable visual artifacts as illustrated in Figure (\ref{fig:bunny_bases_comparision}). These issues can be mitigated by adjusting the support map parameters or using predefined support regions. Additional measures are necessary to preserve the rotational features of the SPLOCS bases. Full implementation of our algorithm can be found on the \href{https://github.com/ShMonem/Snapshots-Reduction-Subspaces-for-Projective-Dynamics}{git repository}.

	\begin{algorithm}[]
		\textbf{Input:} {\color{red}non-rigid-aligned} centered vertex positions snapshots $\mathcal{S} = [q(t_0), q(t_1), \cdots, q(t_T)]$ and mass matrix $M$. \newline
		\textbf{Initialize:} Bases $U = []$. \newline
		
		Mass {\color{red} weight} snapshots: $\mathcal{S} \leftarrow  M^{\frac{1}{2}} \mathcal{S}$\newline
		\For{ k $\leq$ no. Components}{
			Find vertex $\color{red} v$ with {\color{red}largest deformation},\\
			Extract $U_k$: $1^{st}$ PCA mode at $\mathcal{S}(v)$ \\
			Localize $U_k$ and project it out of $\mathcal{S}$\\
			$U \leftarrow  [U, U_k]$ 
		}
		$ U_{pca}  \leftarrow M^{-\frac{1}{2}} U $. \\
		\If{ compute Splocs}{
			\For{$U_k$ in U}{
				$ {\color{red} id_k} = $ vertex with {\color{red}largest deformation}. \\
				Compute support map for $ {\color{red} id_k}$\\
				Sparsify $ {\color{red} U_k}$ \\
			}
			Optimize energy {\color{red} ADMM}\\
			$U_{splocs} \leftarrow  M^{-\frac{1}{2}} U$
		}
		\caption{Volkwein weighted PCA and SPLOCS bases}
		\label{alg: PCA_SPLOCS}
		
	\end{algorithm}
	
	\begin{figure}[h]		
		\begin{tikzpicture}[scale = 0.8]
			\begin{axis}[ylabel=Global step ($\%$), ymin=0, ymax=12e-2, xmin= 0, xmax=220,]
				\addplot + [restrict x to domain=10:200][mark=o, color=teal,very thick,smooth, line width=1pt,solid]
				table [y=podPosGlobal_relative,x=basis,col sep=comma]
				{GraphPath/bunnyTimeMeasuresSum.csv};
				
				\addplot + [restrict x to domain=10:200][mark=o, color=orange,very thick,smooth, line width=1pt,solid]
				table [y=splocsPosGlobal_relative,x=basis,col sep=comma]
				{GraphPath/bunnyTimeMeasuresSum.csv};
				
				\addplot + [restrict x to domain=10:200][mark=o, color=red,very thick,smooth, line width=1pt,solid]
				table [y=lbsPosGlobal_relative,x=basis,col sep=comma]
				{GraphPath/bunnyTimeMeasuresSum.csv};
				
				\legend{PCA, SPLOCS, LBS}
			\end{axis}
		\end{tikzpicture}
		\caption{Comparing relative time required for the global linear system solve at different bases types.}
		\label{fig: GStimeBunny}
	\end{figure}
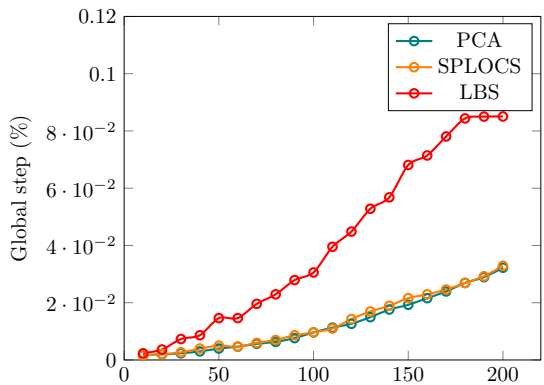

	We discovered that the rigid-pre-alignment of collected snapshots eliminates rotational motion from the information being fed to the computed bases. Hence, we chose to incorporate rotation in the snapshots by avoiding the pre-alignment step. In addition, we employed mass-weighted bases to derive the matrix $U$, as outlined in the Algorithm.


\section{Conclusion}

	We introduce a model reduction technique that provides gain both a speed-up and more realistic simulations. Our snapshots differ from modal analysis because the modes forming the bases are computed from time-varying snapshots, and are not only characterized by the rest state. 
	
\begin{figure}
		\vspace{-1.3cm}
		{\includegraphics[trim={0cm 8cm 0cm 7cm}, clip,width=0.44\textwidth,height=6cm, keepaspectratio]{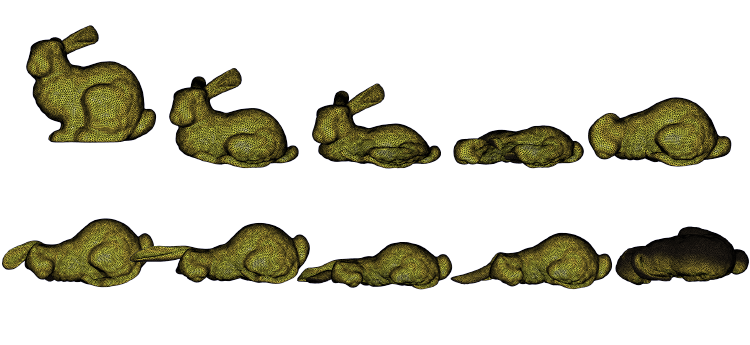}
			\\ \textbf{Full Simulations} 
			\par\vfill
		}
		{\includegraphics[trim={0cm 8cm 0cm 7cm}, clip,width=0.44\textwidth,height=6cm, keepaspectratio]{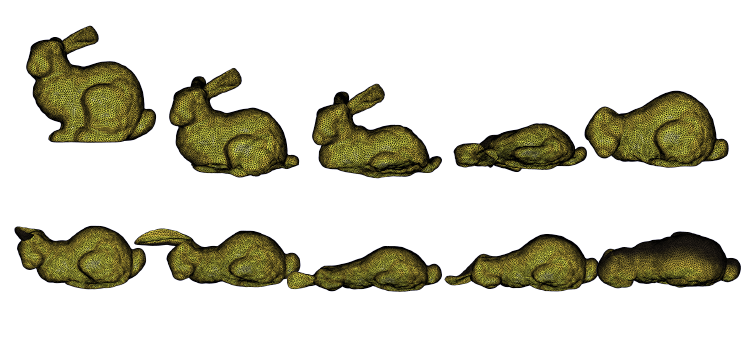} 
			\\ \textbf{PCA}
			\par\vfill
		}
		{\includegraphics[trim={0cm 8cm 0cm 7cm}, clip,width=0.44\textwidth,height=6cm, keepaspectratio]{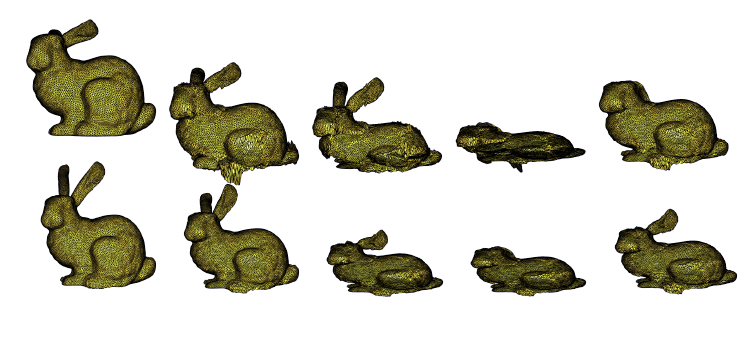}
			\\ \textbf{SPLOCS}
			\par\vfill
		}
		{\includegraphics[trim={0cm 8cm 0cm 7cm}, clip,width=0.44\textwidth,height=6cm, keepaspectratio]{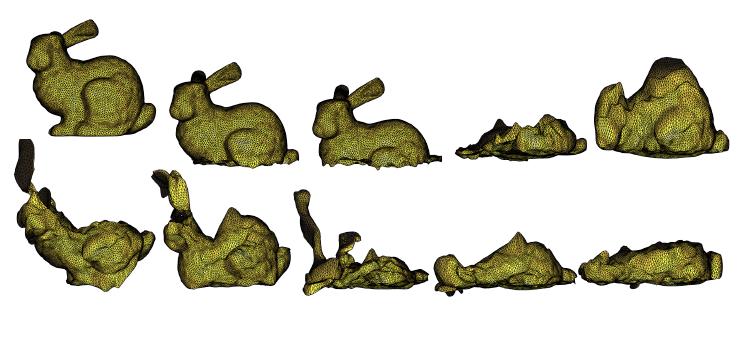}
			\\ \textbf{LBS}
		}  
		\caption{Visual and numerical stability comparison between full model and different reduced simulations, at 200 bases. }
		\label{fig:bunny_bases_comparision}
		
	\end{figure}

	We demonstrated that localized PCA bases from snapshots provide a better representation of global rotations are closest to the original simulations, Figures (\ref{fig:FULL_PCA_SPLOCS_LBS_posOnly}, \ref{fig:bunny_bases_comparision}) and improved computational efficiency for projective dynamics global solve Figure (\ref{fig: GStimeBunny}). Although PCA and SPLOCS bases may appear similar, Figure (\ref{fig:bunny_components_PCA_SPLOCS}), further sparsifying the bases as in SPLOCS can result in the loss of rotational information. 
	
	In summary, our method shows enhanced physical behavior and more efficient global linear solves with improved time margins when using sparse PCA and SLPOCS instead of LBS in projective dynamics for real-time simulations. Only the PCA subspace was capable of representing both rotation and translation efficiently. In conclusion, we computed bases that satisfy
	\begin{itemize}
		\item Visually appealing simulations, realistic and very similar to full order simulations,
		\item Capture time-varying dynamics, particularly potential large deformations,
		\item Faster than state-of-the-art, and the only required extra step is the pre-collected snapshots, by recording the satisfying range of motion from character tests,
		\item Insures high accuracy and numerical stability.
	\end{itemize}
	
	Going forward, we plan to investigate similar bases to further reduce constraint projection.
	
\addcontentsline{toc}{section}{References}
\bibliographystyle{plainurl}
\bibliography{exampleref}

\end{document}